\documentclass{article}

\usepackage{a4wide}
\usepackage{amsmath}
\usepackage{amssymb}
\usepackage{amsthm}
\usepackage{txfonts}
\usepackage{mathtools}
\usepackage{mdwlist}
\usepackage{subcaption}

\newtheorem{remark}{Remark}
\newtheorem{theorem}{Theorem}
\newtheorem{ass}{Assumption}

\newcommand{\pt}{\partial}
\newcommand{\laplace}{\Delta}
\newcommand{\eps}{\varepsilon}
\newcommand{\RR}{\mathbb{R}}
\renewcommand{\phi}{\varphi}
\renewcommand{\rho}{\varrho}
\renewcommand{\theta}{\vartheta}
\newcommand{\scal}[2]{\left(#1,#2\right)}
\newcommand{\dual}[2]{\left\langle#1,#2\right\rangle}
\newcommand{\bigscal}[2]{\bigl(#1,#2\bigr)}
\newcommand{\bigdual}[2]{\bigl\langle#1,#2\bigr\rangle}
\newcommand{\m}[1]{\mathcal{#1}}

\newcommand{\mL}{\mathcal{L}}

\newcommand{\E}{\mathrm{e}}
\newcommand{\D}{\mathrm{  d}}

\newcommand{\abs}[1]{\left|#1\right|}
\newcommand{\norm}[1]{\left\|#1\right\|}
\newcommand{\bignorm}[1]{\big\|#1\big\|}

\title{%
  A unified approach to maximum-norm a posteriori error estimation
  for second-order time discretisations of parabolic equations}

\author{Torsten Lin\ss\thanks{Fakult\"at f\"ur Mathematik und Informatik,
        FernUniversit\"at in Hagen,
        Universit\"atsstra{\ss}e 11,
        58095 Hagen,
        Germany,
        \texttt{[torsten.linss,martin.ossadnik]@fernuni-hagen.de}}
   \and Martin Ossadnik\footnotemark[1]
   \and Goran Radojev\thanks{Department of Mathematics and Computer Science, Faculty of Sciences,
        University of Novi Sad, Trg Dositeja Obradovi\'ca~4, 21000 Novi Sad,
        Serbia,
        \texttt{goran.radojev@dmi.uns.ac.rs}.
        GR acknowledges financial support from DAAD and FernUniversität in Hagen
        through visiting grants.}}

\begin{document}

\maketitle

\begin{abstract}
  A class of linear parabolic equations are considered.
  We derive a common framework for the a posteriori error analysis
  of certain second-order time discretisations combined with finite element
  discretisations in space.
  In particular we study the Crank-Nicolson method, the extrapolated Euler
  method, the backward differentiation formula of order~2 (BDF-2), the
  Lobatto IIIC method and a two-stage SDIRK method.
  We use the idea of elliptic reconstructions and certain bounds
  for the Green’s function of the parabolic operator.
  

\end{abstract}

\section{Introduction}

Given a second-order linear elliptic operator $\m{L}$ in a spatial domain
$\Omega\subset\RR^n$ with Lipschitz boundary, we consider the linear
parabolic equation:
\begin{subequations}\label{problem}
\begin{alignat}{2}
 \m{K}u \coloneqq \pt_t u  + \m{L} u & = f\,,
            &\quad& \text{in} \quad Q\coloneqq\Omega \times (0,T],\\
 \intertext{subject to the initial condition}
   u(x,0) & = u^0(x)\,, && \text{for} \quad x\in\bar{\Omega}, \\
 \intertext{and a homogeneous Dirichlet boundary condition}
   u(x,t) & = 0\,, && \text{for} \quad (x,t)\in \pt\Omega \times [0,T].
\end{alignat}
\end{subequations}
%

There is a vast literature dealing with numerical methods for approximating
parabolic problems.
One classical example is Thom\'ee's monograph \cite{MR2249024} which gives a
comprehensive overview of the FEM for~\eqref{problem} and related equations.
The subject is also elaborated on in various textbooks, e.g.~\cite{MR2362757,
MR3015004}.
The majority of these publications focus on \emph{a priori} error estimation,
i.e. convergence results given in terms of the mesh size an the regularity of
the \emph{exact} solution.

More recently, the derivation of \emph{a posteriori} error bounds has attracted
the attention of many researchers.
This approach yields error bounds in terms of the \emph{computed} numerical
solution, and therefore yield computable bounds on the error.

Most publications to date study the error in $L_2$-norms or in energy norms
naturally induced by the variational formulation of the problem, see
e.g.~\cite{MR695602,MR695603,MR1043607}.
In contrast, our focus is on a posteriori error bounds in the \emph{maximum norm}.
The first such results were probably given in~\cite{MR1335652}, however the
proofs were deferred to a later paper which seems to have never been published.

A key publication is~\cite{MR2034895} by Makridakis and Nochetto who introduced
the concept of \emph{elliptic reconstructions}.
This idea was used in a number of publications to study the FEM combined with
various time discretisations: backward Euler, Crank Nicolson, discontinuous
Galerkin and extrapolation~\cite{MR2519598, MR2629992, MR3056758, AX2208.08153}.
The techniques in those papers are all tailored to the particular method(s)
under consideration and therefore differ to some extend.
Moreover, different -- although related -- stability results have been used.
It is the primary aim of this study to present a \emph{common framework} for the
a posteriori error analysis of second order time discretisations.
Furthermore, all constants appearing in the parabolic error estimator will be
given explicitly.

The paper is organised as follows.
In Section~\ref{sect:weak-discr} we specify our assumptions on the data of
problem~\eqref{problem} and recapitulate certain aspects of the existence theory.
In Section~\ref{sect:ana} we derive our framework for the a posteriori error
analysis of a family of discretisations.
We formulate our assumptions for the analysis in \S\ref{ssect:apost-ell}
(existence of error estimators for the elliptic problems) and
in \S\ref{ssect:green} (bounds on the Green's function of the parabolic problem).
In \S\ref{ssect:reco} the concept of elliptic reconstructions is introduced,
while the main result, Theorem~\ref{theo:theo} is derived in \S\ref{ssect:apost-para}.
In \S\ref{sect:apost-appl} these results are applied to a variety of time
discretisations: backward Euler, Crank-Nicoloson, extrapolation, the
backward-differentiation formula of order $2$, the Lobatto IIIC method and
to a two-stage SDIRK method, with the latter three being novel results.
Finally, numerical results are presented in Section~\ref{sect:numer} to illustrate our
theoretical findings.

\section{Weak formulation and discretisation}
\label{sect:weak-discr}

We shall study~\eqref{problem} in its standard variational form,
cf.~\cite[\S5.1.1]{MR2362757}.
The appropriate Gelfand triple consists of the spaces
\begin{gather*}
  V = H_0^1(\Omega), \quad H= L_2(\Omega) \quad\text{and}\quad
  V^* = H^{-1}(\Omega)\,.
\end{gather*}
Moreover, by \mbox{$a(\cdot,\cdot) \colon V \times V \to \RR$}
we denote the bilinearform associated with the elliptic operator $\m{L}$,
while $\dual{\cdot}{\cdot}\colon V^*\times V \to \RR$ is the duality
pairing and $\scal{\cdot}{\cdot}\colon H\times H\to\RR$ is the scalar
product in~$H$.
Also we shall denote by $\norm{\cdot}_{q,\Omega}$ the standard
norm in $L_q(\Omega)$, \mbox{$q\in[1,\infty]$}.

The solution $u$ of~\eqref{problem} may be considered as a mapping
\mbox{$[0,T] \to V\colon t \mapsto u(t)$}, and we will denote its
(temporal) derivative by $u'$ (and $\pt_t u$).
Let
\begin{gather*}
  W_2^1(0,T;V,H) \coloneqq \left\{ v \in L_2(0,T;V) \colon
         v'\in L_2(0,T;V^*)\right\}\,.
\end{gather*}

Our variational formulation of~\eqref{problem} reads:
Given \mbox{$u^0\in H$} and \mbox{$F\in L_2(0,T;V^*)$},
find \mbox{$u\in W_2^1(0,T;V,H)$} such that
\begin{subequations}\label{weak}
\begin{alignat}{2}
   \frac{\D}{\D t} \bigscal{u(t)}{\chi} + a\bigscal{u(t)}{\chi}
            & = \bigdual{F(t)}{\chi} \quad \forall \chi\in V,
      \ \ t\in(0,T], \\
  \intertext{and}
    u(0)=u^0.
\end{alignat}
\end{subequations}
This problem possesses a unique solution.

In the sequell we shall assume that the source term $F$ has more
regularity and can be represented as \mbox{$\dual{F(t)}{\chi}=\scal{f}{\chi}$},
\mbox{$\forall v\in V$}, with a function \mbox{$f\in L_2(0,T;H)$}

Since we are interested in maximum-norm error estimates we have to make further
assumptions on the data to ensure that the solution can be evaluated pointwise.
To this end, we assume that the intial and boundary data satify the zero-th
order compatibility condition, i.e. $u^0 = 0$ on $\pt\Omega$, and that $u^0$
is Hölder continuous in $\bar\Omega$.
Under standard assumptions on $f$ and $\m{L}$, problem~\eqref{problem} possesses
a unique solution that is continuous on $\bar{Q}$; see~\cite[\S5, Theorem 6.4]{0174.15403}.

Now we turn to discretising~\eqref{weak}.
To this end, let the mesh in time be given by
\begin{gather*}
  \omega_t \colon 0=t_0<t_1<\ldots<t_M=T,
  \intertext{with mesh intervals}
     I_j\coloneqq(t_{j-1},t_j) \quad
  \text{and step sizes} \quad
  \tau_j\coloneqq t_j-t_{j-1},\ j=1,2,\dots,M.
\end{gather*}
For any function \mbox{$v\colon \Omega\times[0,T]\to\RR$} that is continuous
in time on \mbox{$[0,T]$} we set
%
\mbox{$v^j \coloneqq v(\cdot,t_j)$} \mbox{$j=0,1,\dots,M$}.
We also set
\begin{gather*}
  \delta_t v^j \coloneqq \frac{v^j-v^{j-1}}{\tau_j}\,, \ \ j=1,2,\dots,M.
\end{gather*}

Let $V_h$ be a finite dimentional (FE-)subspace of $V$ and
let $a_h\scal{\cdot}{\cdot}$ and $\scal{\cdot}{\cdot}_h$
be approximations of the bilinear form $a\scal{\cdot}{\cdot}$
and of the scalar product $\scal{\cdot}{\cdot}$ in $H$.
These may involve quadrature, for example.

We shall consider an arbitrary time-stepping procedure that generates as
sequence \mbox{$u^j_h\in V_h$}, \mbox{$j=0,1,\dots,M$}, of approximations
to the solution~$u$ of~\eqref{problem} at time $t_j$,
i.e. \mbox{$u^j_h\approx u(t_j)$}.

\section{A posteriori error analysis}
\label{sect:ana}

Our analysis uses three main ingredients:
\begin{itemize*}
  \item a posteriori error bounds for the elliptic problem $\m{L}y=g$,
    see \S\ref{ssect:apost-ell},
  \item bounds for the \textsc{Green}'s function associated with the
    parabolic operator $\m{K}$, see \S\ref{ssect:green} and
  \item the idea of elliptic reconstructions introduced by
    Makridakis and Nochetto~\cite{MR2034895}, see \S\ref{ssect:reco}.
\end{itemize*}
After these concepts have been reviewed, we derive a common framework for
the a posteriori error estimation for second-order time discretisations (and
FEM in space) for~\eqref{problem} in~\S\ref{ssect:apost-para}.

\subsection{A posteriori error estimation for the elliptic problem}
\label{ssect:apost-ell}

Given $g\in H$, consider the elliptic boundary-value problem of
finding $y\in V$ such that
\begin{gather}\label{prob:ell}
   a\scal{y}{\chi} = \scal{g}{\chi}\,, \ \ \forall \chi\in V,
\end{gather}
and its discretisation of finding $y_h\in V_h$ such that
\begin{gather}\label{FEM:ell}
   a_h\scal{y_h}{\chi} = \scal{g}{\chi}_h\,, \ \ \forall \chi\in V_h.
\end{gather}

\begin{ass}\label{ass:ee}
  There exists an a posteriori error estimator $\eta$ for the
  FEM~\eqref{FEM:ell} applied to the elliptic problem~\eqref{prob:ell} with
  \begin{gather*}
    \norm{y_h-y}_{\infty,\Omega} \le \eta\bigl(y_h, g\bigr).
  \end{gather*}
\end{ass}

A few error estimators of this type are available in the literature.
We mention some of them.
\begin{itemize*}
  \item Nochetto et al.~\cite{05068766} study the semilinear problem
        $-\laplace u + g(\cdot,u)=0$ in up to three space dimensions.
        They give a posteriori error bounds for arbitrary order FEM
        on quasiuniform triangulations.
  \item Demlow \& Kopteva~\cite{MR3520007} too consider arbitrary order FEM
        on quasiuniform triangulations, but for the singularly perturbed
        equation $-\eps^2\laplace u + g(\cdot,u)=0$.
        A posteriori error estimates are established that are robust in the
        perturbation parameter.
        Furthermore, in~\cite{MR3419888} for the same problem $P_1$-FEM on
        \emph{anisotropic} meshes are investigated.
  \item In~\cite{MR2334045,MR3232628} arbitrary order FEM for
        the linear problem \mbox{$-\eps^2 u'' + ru = g$} in \mbox{$(0,1)$},
        \mbox{$u(0)=u(1)=0$} are considered.
        In contrast to the afore mentioned contributions all constants
        appearing in the error estimator are given explicitly.
\end{itemize*}

\subsection{Green's functions}
\label{ssect:green}

Let the \textsc{Green}'s function associated with $\m{K}$ and an arbitrary
point \mbox{$x\in\Omega$} be denoted by $\m{G}$,
Then for all \mbox{$\phi\in W^1_2\left(0,T;V,H\right)$}
\begin{gather}\label{green-rep}
  \phi(x,t) = \bigscal{\phi(0)}{\m{G}(t)}
              + \int_0^t \bigdual{\bigl(\m{K}\phi\bigr)(s)}{\m{G}(t-s)} \D s,
    \ \ \text{see~\cite{0144.34903}.}
\end{gather}
The \textsc{Green}'s function
\mbox{$\m{G} \colon \bar{\Omega}\times[0,T]\to\RR$}, \mbox{$t \in (0,T]$}.
solves, for fixed $x$,
\begin{gather*}
  \pt_t\m{G}+\m{L}^*\m{G} = 0, \ \ \text{in} \ \Omega\times\RR^+, \ \
  \m{G}\bigr|_{\pt\Omega} = 0, \ \ \m{G}(0) = \delta_x=\delta(\cdot-x)\,.
\end{gather*}

\begin{ass}\label{ass:green}
  There exist non-negative constants
  $\kappa_0$, $\kappa_1$, $\kappa_1'$ and $\gamma$ such that
  \begin{gather}\label{source:ass}
    \norm{\m{G}(t)}_{1,\Omega} \le \kappa_0\,\E^{-\gamma t} \eqqcolon \phi_0(t),
       \quad
    \norm{\pt_t \m{G}(t)}_{1,\Omega}
       \le \left(\frac{\kappa_1}{t} +\kappa_1'\right)
                  \,\E^{-\gamma t} \eqqcolon \phi_1(t),
  \end{gather}
  for all $x\in\bar\Omega$, $t\in[0,T]$.
\end{ass}

In \S\ref{sect:numer} we will present numerical results for an example
test problem that satisfies these assumptions.
A more detailed discussion of problem classes for which such results are
available is given in~\cite[\S2]{MR3720388}, see also Appendix A in~\cite{MR3056758}.

\subsection{Elliptic reconstruction}
\label{ssect:reco}

Given an approximation $u_h^j\in V_h$ of $u(t_j)$, we define
$\psi^j\in V_h$ by
\begin{gather}\label{reconstr-psi}
  \scal{\psi^j}{\chi}_h = a_h\scal{u_h^j}{\chi}
              - \scal{f^j}{\chi}_h \quad \forall \ \chi\in V_h\,,
      \ \ j= 0,\dots,M.
\end{gather}
This can be written as an ``elliptic'' problem:
\begin{gather}\label{elliptic:discr}
  a_h\scal{u_h^j}{\chi} = \scal{f^j + \psi^j}{\chi}_h \quad
      \forall \ \chi\in V_h, \ \ j=0,\dots,M.
\end{gather}
Next, define $R^j\in H_0^1(\Omega)$ by
\begin{gather}\label{reconstr}
  a\scal{R^j}{\chi} = \scal{f^j+\psi^j}{\chi}
    \quad \forall \ \chi\in H_0^1(\Omega)\,, \ \ j=0,\dots,M,
\end{gather}
or for short:
\begin{gather}\label{reco-LRfpsi}
  \m{L}R^j = f^j+\psi^j\,, \ \ j=0,\dots,M.
\end{gather}
The function $R^j$ is referred to as the elliptic reconstruction of
$u_h^j$, \cite{MR2034895}.


Now, $u_h^j$ can be regarded as the finite-element approximation of $R^j$
obtained by~\eqref{elliptic:discr},
and the error \mbox{$E^j\coloneqq R^j-u_h^j$} can be bounded using the
elliptic estimator from~\S\ref{ssect:apost-ell}:
\begin{gather}\label{reconstr-esti}
  \norm{E^j}_\infty =
  \norm{u_h^j - R^j}_\infty
    \le \eta_\mathrm{ell}^j
    \coloneqq \eta\left(u_h^j, f^j+\psi^j\right)\,,
      \ \ j=0,\dots,M.
\end{gather}
Because of linearity, we have
\begin{gather}\label{reconstr-delta-esti}
  \norm{\delta_t E^j}_\infty =
  \norm{\delta_t\left(u_h - R\right)^j}_\infty
    \le \eta_\mathrm{ell,\delta}^j
    \coloneqq \eta\left(\delta_t u_h^j, \delta_t\left(f+\psi\right)^j\right)\,,
      \ \ j=1,\dots,M.
\end{gather}

\subsection{A posteriori error estimation for the parabolic problem
 -- general framework}
\label{ssect:apost-para}

We are now in a position to derive the main result of the paper.
We like to use the \textsc{Green}'s function representation~\eqref{green-rep}
with $\phi$ replaced by the error \mbox{$u-u_h$}.
First we have to extend the \mbox{$u_h^j$}, \mbox{$j=0,1,\dots,M$}, to a
function defined on all of $[0,T]$.
We use piecewise linear interpolation:
For any function $\phi$ defined on $\omega_t$, $t_j\mapsto \phi^j$,
we define
\begin{gather}\label{hat_notation}
  \hat{\phi}(\cdot,t) \coloneqq \frac{t_{j}-t}{\tau_j}\, \phi^{j-1}
  + \frac{t-t_{j-1}}{\tau_j}\,\phi^{j}
  \quad\text{for}\ \ t\in[t_{j-1},t_j],\quad j=1,\dots, M.
\end{gather}

Eq.~\eqref{green-rep} yields for the error at final time $T$ and for
any $x\in\Omega$:
\begin{gather}\label{rep-error}
  \left(u-u_h^M\right)(x) =
  \left(u-\hat{u}_h\right)(x,T) = \bigscal{u^0-u_h^0}{\m{G}(T)}
              + \int_0^t \bigdual{\bigl(\m{K}(u-\hat{u}_h\bigr)(s)}{\m{G}(t-s)} \D s.
\end{gather}

For the residual of $\hat{u}_h$ in the
differential equation we have the representation
\begin{align*}
  \m{K}(u-\hat{u}_h)(t)
    & = f(t) - \mL\hat{u}_h(t) - \pt_t \hat{u}_h(t) \\
    & = \bigl(f-\hat{f}\bigr)(t)
             + \mL\left(\hat{R} - \hat{u}_h\right)(t)
         - \hat{\psi}(t) - \pt_t \hat{u}_h(t)\,,
\end{align*}
by~\eqref{reco-LRfpsi}.
Substituting this into~\eqref{rep-error} and recalling that \mbox{$E=R-u_h$},
we obtain
\begin{align}\label{error-rep-discr}
  \begin{split}
    \left(u-\hat{u}_h\right)(x,T)
       & = \scal{u^0-u^0_h}{\m{G}(T)}
            + \int_{0}^{T}\scal{\bigl(f-\hat{f}\bigr)(t)}{\m{G}(T-t)} \D t \\
       & \qquad
         + \int_0^T \dual{\mL\hat{E}(t)}{\m{G}(T-t)} \D t
         - \int_0^T \dual{\left(\hat{\psi}+\pt_t\hat{u}_h\right)(t)}{\m{G}(T-t)} \D t\,.
  \end{split}
\end{align}

\begin{theorem}\label{theo:theo}
  Let \mbox{$u^j_h\in V_h$}, \mbox{$j=0,1,\dots,M$}, be an arbitrary sequence
  of approximations to $u(t_j)$.
  Then, for any $K\in\{0,\dots,M-1\}$, one has
  \begin{gather*}
    \norm{u(T)-u_h^M}_{\infty,\Omega}
      \le \eta^{M,K}
      \coloneqq \eta_{\mathrm{init}} + \eta_\mathrm{ell}^{M,K}
        + \eta_f + \eta_{\delta\psi} + \eta_{\Psi},
  \end{gather*}
  where the components of the error estimator $\eta^{M,K}$ are
  \begin{gather*}
    \eta_{\mathrm{init}}\coloneqq \kappa_0 \sigma_0 \norm{u^0-u_h^0}_{\infty,\Omega}, \quad
    \eta_f \coloneqq \sum_{j=1}^{M} \sigma_j \kappa_0
                \int_{I_j} \norm{(f-\hat{f})(s)}_{\infty,\Omega} \D s, \quad
    \eta_{\delta\psi} \coloneqq \sum_{j=1}^{M} \sigma_j \chi_j
                \norm{\delta_t \psi^j}_{\infty,\Omega}\,, \\
    \eta_{\Psi} \coloneqq \kappa_0 \sum_{j=1}^{M} \sigma_j \tau_j \norm{\Psi^j}_{\infty,\Omega}\ \ \
      \text{with}\ \ \
    \Psi^j \coloneqq \frac{\psi^j+\psi^{j-1}}{2} + \delta_t u_h^j\,, \\
   \intertext{and}
    \eta_\mathrm{ell}^{M,K}
        \coloneqq \kappa_0\left(\eta_\mathrm{ell}^M + \sigma_K \eta_\mathrm{ell}^K
         + \sum_{j=K+1}^M \sigma_j \tau_j \eta_{\mathrm{ell},\delta}^j \right)
         + \sum_{j=1}^K 
            \sigma_j \mu_j \max\left\{\eta_\mathrm{ell}^j, \eta_\mathrm{ell}^{j-1}\right\}\,.
  \end{gather*}
  The weights are given by \mbox{$\sigma_j \coloneqq\E^{-\gamma(T-t_j)}$},
  \begin{gather*}
    \mu_j\coloneqq \int_{I_j} \left(\frac{\kappa_1}{T-s} + \kappa_1'\right) \D s \ \ \
    \text{and} \ \ \
    \chi_j \coloneqq \min\left\{\frac{\kappa_0 \tau_j^2}{4},
                             \int_{I_j} \frac{(t_j-s)(s-t_{j-1})}2
                                \left(\frac{\kappa_1}{T-s} + \kappa_1'\right) \D s
               \right\}.
  \end{gather*}
  The elliptic estimators $\eta_\mathrm{ell}^j$ and
  $\eta_{\mathrm{ell},\delta}^j$ have been defined in
  \eqref{reconstr-esti} and~\eqref{reconstr-delta-esti}.
\end{theorem}
\begin{proof}
  We have to bound the right-hand side of~\eqref{error-rep-discr}
  and consider the various terms separately.
  To this end H\"older's inequality and \eqref{source:ass} will be used
  repeatedly.

  \paragraph{(i)}
  First,
  \begin{gather}\label{est-init}
    \abs{\scal{u^0-u_h^0}{\m{G}(T)}}
       \le \kappa_0 \E^{-\gamma T} \eta_{\mathrm{init}}
    \intertext{and}\label{est-F}
    \abs{\int_0^T \scal{(f-\hat{f})(s)}{\m{G}(T-s)}\D s}
       \le \kappa_0 \sum_{j=1}^M \E^{-\gamma(T-t_j)} \eta_f^j\,.
  \end{gather}

  \paragraph{(ii)}
  For the third term on the right-hand side of~\eqref{error-rep-discr},
  we have
  \begin{gather*}
    \int_0^T \scal{\m{L}\bigl(\hat{R}-\hat{u}_h\bigr)(s)}{\m{G}(T-s)}\D s
      = \int_0^T \scal{\pt_t \m{G}_t(T-s)}{\bigl(\hat{R}-\hat{u}_h\bigr)(s)}\D s,
  \end{gather*}
  because $\m{L}^*\m{G} = \pt_t \m{G}$.
  For any $K\in\{0,\dots,M-1\}$, integration by parts on $(t_K,T)$ , gives
  \begin{align*}
    & \int_0^T \scal{\pt_t\m{G}(T-s)}{\bigl(\hat{R}-\hat{u}_h\bigr)(s)}\D s \\
    & \qquad
       = - \scal{\m{G}(0)}{\bigl(R-u_h\bigr)^M}
         + \scal{\m{G}(T-t_K)}{\bigl(R-u_h\bigr)^K}
         + \sum_{j=K+1}^M \int_{I_j}
           \scal{\m{G}(T-s)}{\delta_t\bigl(R-u_h\bigr)^j} \D s \\
    & \qquad\qquad
         + \sum_{j=1}^K \int_{I_j} \scal{\pt_t\m{G}(T-s)}{\bigl(\hat{R}-\hat{u}_h\bigr)(s)}\D s
  \end{align*}
  We apply H\"older's inequality, \eqref{source:ass},
  \eqref{reconstr-esti} and~\eqref{reconstr-delta-esti}
  to obtain
  \begin{gather}\label{est-L(R-u)}
    \begin{split}
    & \abs{\int_0^T \dual{\m{L}\bigl(\hat{R}-\hat{u}_h\bigr)(s)}{\m{G}(T-s)}\D s} \\
    & \qquad
       \le \kappa_0\left(\eta_\mathrm{ell}^M + \E^{-\gamma(T-t_K)} \eta_\mathrm{ell}^K
         + \sum_{j=K+1}^M \E^{-\gamma(T-t_j)} \tau_j \eta_{\mathrm{ell},\delta}^j \right)
         + \sum_{j=1}^K
            \int_{I_j} \phi_1(T-s) \D s \
                   \max\left\{\eta_\mathrm{ell}^j, \eta_\mathrm{ell}^{j-1}\right\}\,.
    \end{split}
  \end{gather}

  \paragraph{(iii)}
  For the last term in~\eqref{error-rep-discr} there holds
  \begin{gather*}
     \left(\hat{\psi}+\pt_t\hat{u}_h\right)(t)
        = \Psi^j + \bigl(t-t_{j-1/2}\bigr) \delta_t \psi^j\,,
                 \ \ t\in\left(t_{j-1},t_j\right), \ \ j=1,\dots,M.
  \end{gather*}
  Using integration by parts, we obtain
  \begin{align*}
    & \int_{t_{j-1}}^{t_j} \dual{\left(\hat{\psi}+\pt_t\hat{u}_h\right)(t)}{\m{G}(T-t)} \D t \\
    & \qquad    =
    \int_{t_{j-1}}^{t_j} \dual{\Psi^j}{\m{G}(T-t)} \D t
       + \int_{t_{j-1}}^{t_j} \bigl(t-t_{j-1/2}\bigr)\dual{\delta_t\psi^j}{\m{G}(T-t)} \D t \\
    & \qquad    =
    \int_{t_{j-1}}^{t_j} \dual{\Psi^j}{\m{G}(T-t)} \D t
       + \frac{1}{2} \int_{t_{j-1}}^{t_j} \bigl(t-t_j\bigr)\bigl(t-t_{j-1}\bigr)
                 \dual{\delta_t\psi^j}{\pt_t \m{G}(T-t)} \D t
  \end{align*}
  Hence
  \begin{align*}
    & \abs{\int_{t_{j-1}}^{t_j}
       \dual{\left(\hat{\psi}+\pt_t\hat{u}_h\right)(t)}
                                     {\m{G}(T-t)} \D t} \\
    & \qquad \le \tau_j \norm{\Psi^j}_{\infty,\Omega} \phi_0(T-t_j)
       + \min\left\{\frac{\tau_j}{4} \phi_0(T-t_j),
                    \frac{1}{2} \int_{t_{j-1}}^{t_j}
                       \bigl(t_j-t\bigr)\bigl(t-t_{j-1}\bigr) \phi_1(T-t_j) \D t \right\}
         \cdot \norm{\delta_t\psi^j}_{\infty,\Omega}\,.
  \end{align*}
  This gives
  \begin{gather}\label{est-Psi}
    \abs{\int_{t_{j-1}}^{t_j} \dual{\left(\hat{\psi}+\pt_t\hat{u}_h\right)(t)}
                                     {\m{G}(T-t)} \D t}
    \le \kappa_0 \tau_j \sigma_j \norm{\Psi^j}_{\infty,\Omega}
       + \chi_j \norm{\delta_t\psi^j}_{\infty,\Omega}
  \end{gather}

  Finally, applying~\eqref{est-init}--\eqref{est-Psi}
  to~\eqref{error-rep-discr} completes the proof.
\end{proof}

\begin{remark}\label{rem:est}
  \emph{(i)}
  In general, the supremum norm involved in $\eta_\mathrm{init}$ can not be
  determined exactly, but needs to be approximated.
  For example, one can use a mesh that is finer than the finite-element mesh.

  \emph{(ii)}
  The integral in $\eta_f^j$ needs to be approximated.
  One possibility is  Simpson's rule, which is of order $4$ and gives
  \begin{gather*}
    \int_{I_j} \norm{\bigl(f-\hat{f}\bigr)(s)}_{\infty,\Omega} \D s
       \approx \frac{\tau_j}{3}
               \norm{f^j - 2 f^{j-1/2} + f^{j-1}}_{\infty,\Omega}
       \approx \frac{\tau_j^3}{12} \norm{\pt_t^2 f(t_{j-1/2})}_{\infty,\Omega}
       \,.
  \end{gather*}
  Here too, the supremum norm needs to be approximated.
\end{remark}


\section{Application to various time discretisations}
\label{sect:apost-appl}

The framework derived in the preceding section has not made any use of a
particular time discretisation.
Theorem~\ref{theo:theo} does not discriminate between them.
The differences become obvious when analysing~\eqref{reconstr-psi},
the definition of the $\psi^j$.
This also reveals alternative formulae for computing the $\psi^j$,
\mbox{$j>0$}, that --
unlike~\eqref{reconstr-psi} -- do not require to invert the mass matrix.

\subsection{The backward Euler method}

This methods reads as follows: Given an approximation $u_h^0$ of the initial
data, find \mbox{$u^j_h \in V_h$}, \mbox{$j = 1,\dots,M$}, such that
\begin{gather}\label{euler-FE}
   \scal{\frac{u_h^j-u_h^{j-1}}{\tau_j}}{\chi}_h + a_h \scal{u_h^j}{\chi}
      = \scal{f^j}{\chi}_h \quad \forall \ \chi \in V_h.
\end{gather}
Comparing this equation with~\eqref{reconstr-psi}, we see that
\mbox{$\psi^j = - \delta_t u_h^j$}, \mbox{$j=1,\dots,M$}, and
\begin{gather*}
  \bigl(\hat{\psi} + \pt_t \hat{u}_h\bigr)(t)
    = \frac{\psi^{j-1}-\psi^j}{2} + \bigl(t-t_{j-1/2}\bigr) \delta_t \psi^j
    = -\frac{\tau_j}{2}\delta_t \psi^j
               + \bigl(t-t_{j-1/2}\bigr) \delta_t \psi^j\,,\ \ j=1,\dots,M.
\end{gather*}
Thus, $\Psi^j = -\frac{\tau_j}{2}\delta_t \psi^j$, \mbox{$j=1,\dots,M$}.

\subsection{The Crank-Nicolson method}

Given $u_h^0\approx u^0$, find \mbox{$u^j_h \in V_h$}, \mbox{$j = 1,\dots,M$},
such that
\begin{gather}\label{CN-FE}
   \scal{\frac{u_h^j-u_h^{j-1}}{\tau_j}}{\chi}_h\
      + a_h \scal{\frac{u_h^j+u_h^{j-1}}{2}}{\chi}
      = \scal{\frac{f^j+f^{j-1}}{2}}{\chi}_h \quad \forall \ \chi \in V_h.
\end{gather}
This method may also be viewed as a Runge-Kutta-Lobatto-IIIA
method~\cite{MR239762}. It is $A$-stable, but not $L$-stable.

Comparison with~\eqref{reconstr-psi} shows that
\mbox{$\bigl(\psi^j+\psi^{j-1}\bigr)/2 = - \delta_t u_h^j$}.
Therefore,
\begin{gather*}
  \bigl(\hat{\psi} + \pt_t \hat{u}_h\bigr)(t)
    = \bigl(t-t_{j-1/2}\bigr) \delta_t \psi^j\,, \ \ j=1,\dots,M.
\end{gather*}
In particular, $\Psi^j = 0$, \mbox{$j=1,\dots,M$},
and the term $\eta_\Psi^j$ disappears from the error estimator.

\subsection{Extrapolation based on the backward Euler method}

Starting from $u_h^0\approx u^0$, three sequences of approximations
are generated as follows.

\begin{subequations}\label{extra-FE}
  \noindent
  \textbf{One-step Euler:}
    Set $v_h^0=u_h^0$ and find $v^j_h \in V_h$, $j=1,\dots,M$, such that
    \begin{align}\label{extra-FE-one}
       \scal{\frac{v_h^j-v_h^{j-1}}{\tau_j}}{\chi}_h + a_h \scal{v^j_h}{\chi}
          & = \scal{f^j}{\chi}_h \quad \forall \ \chi \in V_h.
  \intertext{\textbf{Two-step Euler:}
    Set $w_h^0=u_h^0$ and find $w_h^{j-1/2},w^j_h \in V_h$, $j=1,\dots,M$,
    such that}
       \label{extra-FE-two-1}
       \scal{\frac{w_h^{j-1/2} - w_h^{j-1}}{\tau_j/2}}{\chi}_h
          + a_h \scal{w_h^{j-1/2}}{\chi}
          & = \scal{f^{j-1/2}}{\chi}_h \quad \forall \ \chi \in V_h, \\
       \label{extra-FE-two-2}
       \scal{\frac{w_h^j - w_h^{j-1/2}}{\tau_j/2}}{\chi}_h
          + a_h \left( w_h^j,\chi\right)
          & = \scal{f^j}{\chi}_h \quad \forall \ \chi \in V_h.
    \end{align}
  \textbf{Extrapolation:} Set
    \begin{gather}
       u_h^j \coloneqq 2 w_h^j - v_h^j, \quad j=1,\dots,M.
    \end{gather}
\end{subequations}
For this method eqs.~\eqref{reconstr-psi}, \eqref{extra-FE-one}
and~\eqref{extra-FE-two-2} yield
\begin{gather*}
  \psi^j = - 4 \frac{w_h^j-w_h^{j-1/2}}{\tau_j}
      + \frac{v_h^j-v_h^{j-1}}{\tau_j}\,, \ \ j=1,\dots,M.
\end{gather*}

\subsection{The backward-differentiation formula (BDF-2)}
Let us briefly recall the constuction of the BDF-2 method.
Given approximations $U^{j-2}$ and $U^{j-1}$ of $u(t_{j-2})$ and $u(t_{j-1})$,
\mbox{$j\ge 2$}, we seak an approximation $U^j$ of $u(t_j)$ as the solution of
\begin{gather*}
  U_2'(t_j) + \m{L} U_2(t_j) = f^j,
\end{gather*}
where $U_2$ is the uniquely defined quadratic interpolation polynomial with
$U_2(t_\iota)=U^\iota$ for $\iota\in\{j-2,j-1,j\}$.
This idea yields the time-stepping procedure
\begin{gather*}
  D_t U^j + \m{L} U^j = f^j, \ \ j=2,\dots,M,
\end{gather*}
where
\begin{gather*}
  D_t v^j \coloneqq \alpha_j \delta_t v^j + \beta_j \delta_t v^{j-1}, \quad
   \alpha_j = \frac{2\tau_j+\tau_{j-1}}{\tau_j+\tau_{j-1}}, \quad
   \beta_j = - \frac{\tau_j}{\tau_j+\tau_{j-1}}\,, \quad j=2,\dots,M.
\end{gather*}
The difference quotient $D_t$ can also be represented as
\begin{gather*}
   D_t v^j = \delta_t v^j + \tau_j \delta^2_t v^j \,,
     \quad \delta^2_t v^j \coloneqq \frac{\delta_t v^j - \delta_t v^{j-1}}{\tau_j+\tau_{j-1}}\,,
       \ \ j=2,\dots,M.
\end{gather*}

The BDF-2 time-stepping procedure requires two starting values.
One is given naturally by the initial condition.
The other one is obtained by applying one step of the backward
Euler method on the first time interval.
The local error of that method is $2$ thus matching the formal
order of the BDF-2 method.

We formally set \mbox{$\delta_t^2 v^1 \coloneqq 0$}.
Then the BDF-2 discretisation reads as follows:
seek \mbox{$u^j_h \in V_h$}, \mbox{$j = 1,\dots,M$}, such that
\begin{gather}\label{bdf-FE}
  \scal{D_t u^j_h}{\chi}_h
     + a_h\scal{u^j_h}{\chi} = \scal{f^j}{\chi}_h \quad \forall
         \ \chi \in V_h.
\end{gather}

For this method, eqs.~\eqref{reconstr-psi} and \eqref{bdf-FE} yield
$\psi^j = - \delta_t u_h^j - \tau_j\delta_t^2 u_h^j$ and
\begin{gather*}
  \bigl(\hat{\psi} + \pt_t \hat{u}_h\bigr)(t)
    = -\frac{\tau_j}{2}\delta_t \psi^j - \tau_j \delta_t^2 u_h^j
               + \bigl(t-t_{j-1/2}\bigr) \delta_t \psi^j\,,\ \ j=1,\dots,M,
\end{gather*}
i.e., $\Psi^j = -\frac{\tau_j}{2}\delta_t \psi^j - \tau_j \delta_t^2 u_h^j$,
$j=1,\dots,M$.

\subsection{The two-stage Lobatto-IIIC method}

This Runge-Kutta method --- proposed in~\cite{MR295582} --- is given by the
Butcher table
\begin{center}
  \begin{tabular}{c|cc}
    $0$ & $1/2$ & $-1/2$ \\
    $1$ & $1/2$ & $1/2$ \\\hline
        & $1/2$ & $1/2$ 
  \end{tabular}
\end{center}
In contrast to the Crank-Nicolson method it is both $A$- and $L$-stable.
It can be formulated as follows.
Given $u_h^0\approx u^0$, find \mbox{$v^j_h,u^j_h \in V_h$}, \mbox{$j = 1,\dots,M$},
such that
\begin{subequations}\label{IIIC-FE}
\begin{align}\label{IIIC-FE-one}
  \scal{\frac{v_h^j-u_h^{j-1}}{\tau_j}}{\chi}_h + \frac{1}{2} a_h \scal{v^j_h-u_h^j}{\chi}
          & = \frac{1}{2}\scal{f^{j-1}-f^j}{\chi}_h\quad \forall \ \chi \in V_h, \\
  \scal{\frac{u_h^j-u_h^{j-1}}{\tau_j}}{\chi}_h + \frac{1}{2} a_h \scal{v^j_h+u_h^j}{\chi}
          & = \frac{1}{2}\scal{f^{j-1}+f^j}{\chi}_h \quad \forall \ \chi \in V_h.
\end{align}
\end{subequations}
Subtracting these two equations, one gets
\begin{gather*}
  \scal{\frac{u_h^j-v_h^{j-1}}{\tau_j}}{\chi}_h + a_h \scal{u_h^j}{\chi}
          = \scal{f^j}{\chi}_h \quad \forall \ \chi \in V_h.
\end{gather*}
Comparison with~\eqref{reconstr-psi} gives
\begin{gather*}
  \psi^j = \frac{v_h^j-u_h^j}{\tau_j}\,, \ \ \ j=1,\dots,M.
\end{gather*}

\subsection{A two-stage SDIRK method}

The singly diagonally implicit Runge-Kutta method given by the
Butcher table
\begin{center}
  \begin{tabular}{c|cc}
    $\gamma$ & $\gamma$ & $0$ \\
    $1-\gamma$ & $1-2\gamma$ & $\gamma$ \\\hline
        & $1/2$ & $1/2$ 
  \end{tabular}
  \qquad with \quad $\displaystyle \gamma=\frac{2-\sqrt{2}}{2}$
\end{center}
is also both $A$- and $L$-stable; see~\cite{MR569982}.
It can be formulated as follows.
Given $u_h^0\approx u^0$, find \mbox{$k^j_1,k^j_2 \in V_h$}, \mbox{$j = 1,\dots,M$},
such that
\begin{subequations}\label{SDIRK-FE}
\begin{align}\label{SDIRK-FE-one}
  \scal{k^j_1}{\chi}_h + \gamma \tau_j a_h \scal{k^j_1}{\chi}
          & = \scal{f^{j-1+\gamma}}{\chi}_h - a_h\scal{u^{j-1}}{\chi}\,, \quad \forall \ \chi \in V_h, \\
   \label{SDIRK-FE-two}
  \scal{k^j_2}{\chi}_h + \gamma \tau_j a_h \scal{k^j_2}{\chi}
          & = \scal{f^{j-\gamma}}{\chi}_h - a_h\scal{u^{j-1}}{\chi}
                 - (1-2\gamma) \tau_j a_h\scal{k^j_1}{\chi}
                      \,, \quad \forall \ \chi \in V_h, \\
  \intertext{and set}
    u^j & = u^{j-1} + \tau_j \frac{k^j_1+k^j_2}{2}\,.
\end{align}
\end{subequations}
Here $f^{j-1+\kappa} \coloneqq f(t_{j-1}+\kappa \tau_j)$ for any $\kappa\in[0,1]$.

The procedure defined by~\eqref{SDIRK-FE} does not involve any evaluation of
$f$ at mesh points $t_j$.
Therefore, in contrast to the methods considered so far, the $\psi^j$ \emph{must}
be computed using~\eqref{reconstr-psi}.

However, if in~\eqref{SDIRK-FE} the function $f$ is replaced by its piecewise
interpolant $\hat{f}$, then one has
\begin{gather*}
  \psi^j = \frac{1-2\gamma}{2\gamma} k^j_1 - \frac{1}{2\gamma} k^j_2
         = \frac{k^j_1-k^j_2}{2\gamma} - k^j_1 \,.
\end{gather*}
To verify this, multiply~\eqref{SDIRK-FE-one} by $(1-1/2\gamma)$,
\eqref{SDIRK-FE-two} by $1/2\gamma$ and sum both equations.

\begin{remark}
  In the above presentation, we have implicitly assumed that
  \eqref{euler-FE}-\eqref{SDIRK-FE} all possess unique solutions.
  This is guarenteed when, for example, the bilinear form $a_h(\cdot,\cdot)$
  is coercive and bounded, and when the $\scal{\cdot}{\cdot}_h$ is a scalar
  product on $V_h$.
  The latter is the case when \mbox{$\scal{\zeta}{\chi}_h=\scal{\zeta}{\chi}$}
  $\forall \zeta,\chi\in V_h$ or when lumping is used.
  These are two standard choices in FE-discretisations of parabolic PDEs.
\end{remark}

\section{A numerical example}
\label{sect:numer}

Consider the following reaction-diffusion equation
\begin{subequations}\label{testproblem}
\begin{alignat}{2}
 \pt_t u  - u_{xx} + (5x+6) u & = \E^{-4t} + \cos\bigl(\pi(x+t)^2\bigr)\,,
            &\quad& \text{in} \quad (-1,1) \times (0,1],\\
 \intertext{subject to the initial condition}
   u(x,0) & = \sin \frac{\pi(1+x)}{2}\,, && \text{for} \quad x\in[-1,1], \\
 \intertext{and the Dirichlet boundary condition}
   u(x,t) & = 0\,, && \text{for} \quad (x,t)\in \{-1,1\} \times [0,1].
\end{alignat}
\end{subequations}
The Green's function for this problem satisfies
\begin{gather}\label{test-green}
  \norm{\m{G}(t)}_{1,\Omega} \le \E^{-t/2},
    \quad
  \norm{\pt_t\m{G}(t)}_{1,\Omega}
            \le \frac{3}{2^{3/2}} \frac{\E^{-t/2}}{t} \,, \ \
    \text{see~\cite[\S12]{MR3056758}}.
\end{gather}
The exact solution to this problem is unknown.
To compute a reference solution, we use a spectral method in space combined
with the dG(2) method in time which is of order $5$. This gives an approximation
that is accurate close to machine precision.

In our experiments, we use the spatial discretisation by $P_1$-FEM analysed
in~\cite{MR2334045} and the a posteriori estimator derived therein.
That method is of order $2$.
Because our time discretisations (except for the Euler method) are also
of second order, we couple spatial and temporal mesh sizes by \mbox{$h=\tau$}.

The maximum norm of the error needs to be approximated.
We do so by evaluating the error at seven evenly distributed points in each
mesh interval:
\begin{gather*}
  \bignorm{u(T)-U^M}_{\infty,\Omega} \approx
  e_M\coloneqq \max_{i=1,\dots,N} \, \max_{r=0,\dots,7}
       \abs{\left(u(T)-U^M\right)\left(x_{i-1} + rh/7\right)}\,,
\end{gather*}
where the $x_i$, $i=0,\dots,N$, are the nodes of the uniform spatial mesh,
and $h$ its mesh size.

Tables~\ref{tab:all}(a-f) display the results of our test
computations.
The first column in each table contains the number of mesh intervals $M$
(with \mbox{$h=\tau=1/M$}),
followed by the errors $e_M$ at final time, the experimental order of
convergence $p_M$, the error estimator $\eta^{M,0}$ and finally
the efficiency $\chi_M$:
\begin{gather*}
  p_M\coloneqq \frac{\ln e_{M/2} - \ln e_M}{\ln 2} \quad \text{and} \quad
  \chi_M \coloneqq \frac{e_M}{\eta^{M,0}}.
\end{gather*}
\begin{table}[h]
\begin{subtable}{.5\textwidth}\centering
\begin{tabular}{cccccccccc}
$M$ & $e_M$ & $p_M$ & $\eta^{M,0}$ & $\chi_M$ \\\hline
    64 & 5.977e-4 & 1.39 & 5.528e-2 & 1/92 \\
   128 & 2.512e-4 & 1.25 & 2.357e-2 & 1/94 \\
   256 & 1.137e-4 & 1.14 & 1.073e-2 & 1/94 \\
   512 & 5.387e-5 & 1.08 & 5.095e-3 & 1/95 \\
  1024 & 2.619e-5 & 1.04 & 2.479e-3 & 1/95 \\
  2048 & 1.291e-5 & 1.02 & 1.223e-3 & 1/95 \\
  4096 & 6.409e-6 & 1.01 & 6.071e-4 & 1/95 \\
  8192 & 3.193e-6 & 1.01 & 3.025e-4 & 1/95 \\
 16384 & 1.594e-6 & 1.00 & 1.510e-4 & 1/95 \\
\hline
\end{tabular}
\caption{Euler}
\end{subtable}
\begin{subtable}{.5\textwidth}\centering
\begin{tabular}{cccccccccc}
$M$ & $e_M$ & $p_M$ & $\eta^{M,0}$ & $\chi_M$ \\\hline
    64 & 2.006e-4 & 1.93 & 1.680e-2 & 1/84 \\
   128 & 5.068e-5 & 1.98 & 4.301e-3 & 1/85 \\
   256 & 1.269e-5 & 2.00 & 1.088e-3 & 1/86 \\
   512 & 3.174e-6 & 2.00 & 2.736e-4 & 1/86 \\
  1024 & 7.935e-7 & 2.00 & 6.863e-5 & 1/86 \\
  2048 & 1.984e-7 & 2.00 & 1.720e-5 & 1/87 \\
  4096 & 4.959e-8 & 2.00 & 4.307e-6 & 1/87 \\
  8192 & 1.240e-8 & 2.00 & 1.078e-6 & 1/87 \\
 16384 & 3.093e-9 & 2.00 & 2.700e-7 & 1/87 \\
\hline
\end{tabular}
\caption{Crank-Nicolson}
\end{subtable}
\begin{subtable}{.5\textwidth}\centering
\begin{tabular}{cccccccccc}
$M$ & $e_M$ & $p_M$ & $\eta^{M,0}$ & $\chi_M$ \\\hline
    64 & 1.986e-4 & 1.92 & 1.873e-2 &  1/94 \\
   128 & 5.024e-5 & 1.98 & 4.849e-3 &  1/97 \\
   256 & 1.259e-5 & 2.00 & 1.240e-3 &  1/99 \\
   512 & 3.148e-6 & 2.00 & 3.155e-4 & 1/100 \\
  1024 & 7.871e-7 & 2.00 & 8.002e-5 & 1/102 \\
  2048 & 1.968e-7 & 2.00 & 2.026e-5 & 1/103 \\
  4096 & 4.919e-8 & 2.00 & 5.126e-6 & 1/104 \\
  8192 & 1.231e-8 & 2.00 & 1.296e-6 & 1/105 \\
 16384 & 3.033e-9 & 2.02 & 3.276e-7 & 1/108 \\
\hline
\end{tabular}
\caption{Extrapolated Euler}
\end{subtable}
\begin{subtable}{.5\textwidth}\centering
\begin{tabular}{cccccccccc}
$M$ & $e_M$ & $p_M$ & $\eta^{M,0}$ & $\chi_M$ \\\hline
    64 & 2.092e-4 & 1.94 & 2.495e-2 & 1/119 \\
   128 & 5.261e-5 & 1.99 & 6.621e-3 & 1/126 \\
   256 & 1.314e-5 & 2.00 & 1.719e-3 & 1/131 \\
   512 & 3.285e-6 & 2.00 & 4.415e-4 & 1/134 \\
  1024 & 8.209e-7 & 2.00 & 1.126e-4 & 1/137 \\
  2048 & 2.052e-7 & 2.00 & 2.862e-5 & 1/139 \\
  4096 & 5.129e-8 & 2.00 & 7.256e-6 & 1/141 \\
  8192 & 1.282e-8 & 2.00 & 1.837e-6 & 1/143 \\
 16384 & 3.205e-9 & 2.00 & 4.648e-7 & 1/145 \\
\hline
\end{tabular}
\caption{BDF-2}
\end{subtable}

\begin{subtable}{.5\textwidth}\centering
\begin{tabular}{cccccccccc}
$M$ & $e_M$ & $p_M$ & $\eta^{M,0}$ & $\chi_M$ \\\hline
    64 & 2.426e-4 & 1.85 & 2.120e-2 & 1/87 \\
   128 & 6.392e-5 & 1.92 & 5.612e-3 & 1/88 \\
   256 & 1.649e-5 & 1.95 & 1.462e-3 & 1/89 \\
   512 & 4.201e-6 & 1.97 & 3.773e-4 & 1/90 \\
  1024 & 1.061e-6 & 1.98 & 9.669e-5 & 1/91 \\
  2048 & 2.669e-7 & 1.99 & 2.467e-5 & 1/92 \\
  4096 & 6.692e-8 & 2.00 & 6.278e-6 & 1/94 \\
  8192 & 1.676e-8 & 2.00 & 1.595e-6 & 1/95 \\
 16384 & 4.185e-9 & 2.00 & 4.044e-7 & 1/97 \\
\hline
\end{tabular}
\caption{Lobatto IIIC}
\end{subtable}
\begin{subtable}{.5\textwidth}\centering
\begin{tabular}{cccccccccc}
$M$ & $e_M$ & $p_M$ & $\eta^{M,0}$ & $\chi_M$ \\\hline
    64 & 2.112e-4 & 1.92 & 1.924e-2 & 1/91 \\
   128 & 5.354e-5 & 1.98 & 4.949e-3 & 1/92 \\
   256 & 1.344e-5 & 1.99 & 1.257e-3 & 1/94 \\
   512 & 3.363e-6 & 2.00 & 3.175e-4 & 1/94 \\
  1024 & 8.412e-7 & 2.00 & 7.996e-5 & 1/95 \\
  2048 & 2.103e-7 & 2.00 & 2.011e-5 & 1/96 \\
  4096 & 5.259e-8 & 2.00 & 5.053e-6 & 1/96 \\
  8192 & 1.315e-8 & 2.00 & 1.269e-6 & 1/97 \\
 16384 & 3.287e-9 & 2.00 & 3.190e-7 & 1/97 \\
\hline
\end{tabular}
\caption{SDIRK}
\end{subtable}
\caption{Error, estimator and efficiency, test problem~\eqref{testproblem}\label{tab:all}}
\end{table}
The methods converge with the expected orders of $1$ (Euler) and $2$ (all others).
The error estimators are upper bounds on the actual errors and both correlate.
Depending on the particular method the errors are overestimated by a factor
of~$100$ to $150$.
Part of this overestimation can be attributed to the bounds
in~\eqref{test-green} not being sharp.
Smaller constants will automatically result in more efficient error
estimators.
The same applies to the constants in the elliptic error estimator $\eta$
in~\S3.1 which is one building block of our Theorem~\ref{theo:theo}.

It is worth studying the various components of the error estimator.
Table~\ref{tab:bdf2-full} displays the results for the BDF-2 discretisation.
While for all four methods the terms $\eta_\mathrm{init}$ and $\eta_f$ are
obviously identical, it is worth noting that the terms $\eta_\mathrm{ell}^{M,0}$
and $\eta_{\delta\psi}$ are also very similar (up to the $4^\mathrm{th}$ digit).
Only the term $\eta_\Psi$ varies signifantly between the methods.
For the Crank-Nicolson method it vanishes identically.
For the Euler method it behaves like $\m{O}(1/M)$, $M\to\infty$.

\begin{table}[!h]
\begin{center}
\centerline{
\begin{tabular}{cccccccccc}
$M$ & $\eta_{\text{init}}$ & $\eta_f$ & $\eta_\text{ell}^{M,0}$ & $\eta_{\Psi}$ & $\eta_{\delta\psi}$ \\\hline
    64 & 1.789e-4 & 1.392e-3 & 1.496e-2 & 8.159e-3 (1.71) & 2.677e-4 (1.89) \\
   128 & 4.474e-5 & 3.474e-4 & 3.836e-3 & 2.321e-3 (1.81) & 7.175e-5 (1.90) \\
   256 & 1.118e-5 & 8.677e-5 & 9.704e-4 & 6.313e-4 (1.88) & 1.913e-5 (1.91) \\
   512 & 2.796e-6 & 2.168e-5 & 2.440e-4 & 1.680e-4 (1.91) & 5.082e-6 (1.91) \\
  1024 & 6.990e-7 & 5.419e-6 & 6.117e-5 & 4.400e-5 (1.93) & 1.345e-6 (1.92) \\
  2048 & 1.748e-7 & 1.355e-6 & 1.531e-5 & 1.142e-5 (1.95) & 3.550e-7 (1.92) \\
  4096 & 4.369e-8 & 3.386e-7 & 3.831e-6 & 2.949e-6 (1.95) & 9.343e-8 (1.93) \\
  8192 & 1.092e-8 & 8.466e-8 & 9.581e-7 & 7.590e-7 (1.96) & 2.453e-8 (1.93) \\
 16384 & 2.731e-9 & 2.116e-8 & 2.396e-7 & 1.949e-7 (1.96) & 6.424e-9 (1.93) \\
 \hline
\end{tabular}
}
\caption{BDF-2: Composition of the error estimator, test problem~\eqref{testproblem}\label{tab:bdf2-full}}
\end{center}
\end{table}

In our test computations, the term $\eta_\mathrm{ell}^{M,0}$ dominates the error
estimator.
This term contains the contributions from the elliptic error estimator.
Most notably the terms $\eta_\Psi$ and $\eta_{\delta\psi}$ do not behave like
$\m{O}(1/M^2)$. There seems to be some kind of logarithmic dependence on the
step size.
Such dependencies are known from other a posteriori error estimates, see
e.g.~\cite{05068766}.

A Matlab/Octave program that reproduces the tables of this paper can be found at
GitHub: \textrm{https://github.com/TorstenLinss/LOR2023}.

\def\cprime{$'$}

\end{document}